\documentclass[a4paper,10pt]{article}

\usepackage{euscript}
\usepackage{amsfonts}
\usepackage{amsmath,amscd, amsthm, amssymb}
\usepackage[english]{babel}
\usepackage[latin1]{inputenc}
\usepackage{graphicx}
\usepackage[arrow, matrix, curve]{xy}
\usepackage{dsfont}\let\mathbb\mathds
\usepackage{rotating}

\newtheorem{thm}{Theorem}

\newtheorem{Def}[thm]{Definition}

\def\CC{{\mathbb{C}}}

\def\RR{{\mathbb{R}}}
\def\ZZ{{\mathbb{Z}}}

\def\PP{{\mathbb{P}}}

\def\cartesien{%
    \ar@{-}[]+R+<6pt,-1pt>;[]+RD+<6pt,-6pt>%
    \ar@{-}[]+D+<1pt,-6pt>;[]+RD+<6pt,-6pt>%
  }
\newcommand{\hook}[1][r]
   {\ar@{}[#1] |*[o][F]{\hbox{%
         \vrule width 1.5mm height 0pt depth 0pt%
         \vrule width 0pt height .75mm depth .75mm%
         }}
     \ar@{^{(}->}[#1]}
\newcommand{\hookl}[1][r]
   {\ar@{}[#1] |*[o][F]{\hbox{%
         \vrule width 1.5mm height 0pt depth 0pt%
         \vrule width 0pt height .75mm depth .75mm%
         }}
     \ar@{_{(}->}[#1]}     
\newcommand{\demde}[1]{\begin{proof} de #1}
\newcommand{\dem}{\begin{proof}}
\newcommand{\cqfd}{\end{proof}}
  \newcommand{\UN}[4][r]{%
    \ar@/^1pc/[#1]^{#2}_*=<0.3pt>{}="HAUT"
    \ar@/_1pc/[#1]_{#3}^*=<0.3pt>{}="BAS"
    \ar @{=>} "HAUT";"BAS" ^{#4}
  }
\newcommand{\DEUX}[6][r]{
    \ar@/^2pc/[#1]^{#2}_*=<0.3pt>{}="HAUT"
    \ar@{}    [#1]     ^*=<0.3pt>{}="MILIEUHAUT"
                       _*=<0.3pt>{}="MILIEUBAS"
    \ar[#1]_(0.3){#3}                  
    \ar@/_2pc/[#1]_{#4}^*=<0.3pt>{}="BAS"
    \ar @{=>} "HAUT";"MILIEUHAUT" ^{#5}
    \ar @{=>} "MILIEUBAS";"BAS" ^{#6}
  }   
 \newcommand{\eq}[1][r]
   {\ar@<-3pt>@{-}[#1]
    \ar@<-1pt>@{}[#1]|<{}="gauche"
    \ar@<+0pt>@{}[#1]|-{}="milieu"
    \ar@<+1pt>@{}[#1]|>{}="droite"
    \ar@/^2pt/@{-}"gauche";"milieu"
    \ar@/_2pt/@{-}"milieu";"droite"}
 \newcommand{\incl}[1][r]
  {\ar@<-0.2pc>@{^(-}[#1] \ar@<+0.2pc>@{-}[#1]}

\newcommand{\Bc}{\mathcal{B}}
\newcommand{\Cc}{\mathcal{C}}

\newcommand{\Ic}{\mathcal{I}}

\newcommand{\Pc}{\mathcal{P}}

\newcommand{\CCC}{\mathfrak{C}}

\title{Faisceaux pervers sur les vari\'et\'es toriques lisses}
\author{Delphine Dupont}
\begin{document}
\maketitle
\begin{abstract}
Soit $X$ une vari\'et\'e torique lisse stratifi\'ee par l'action du tore. \`A partir de l'\'eventail associ\'ee \`a $X$ nous d\'efinissons une cat\'egorie de repr\'esentations de carquois \'equivalente \`a la cat\'egorie des faisceaux pervers sur $X$ relativement \`a cette stratification.
\end{abstract}
~\\
\section{Introduction}
Depuis la d\'efinition dans \cite{BBD} de la cat\'egorie $\Pc erv_X$ des faisceaux pervers sur un espace topologique $X$, de nombreuses m\'ethodes ont \'et\'e d\'evelopp\'ees pour en donner dans le cas o\`u la stratification est fix\'ee une description \'el\'ementaire. Plusieurs descriptions par des cat\'egories de repr\'esentations de carquois se font sur des espaces topologiques particuliers munis d'une stratification fix\'ee. Rappelons qu'un carquois est un graphe orient\'e et qu'une repr\'esentation de ce carquois est un foncteur de ce graphe dans la cat\'egorie des espaces vectoriels de dimension finie. Un premier exemple de cette approche est donn\'e par A. Galligo, M. Granger et Ph. Maisonobe qui, dans \cite{GGM}, d\'emontrent l'\'equivalence de la cat\'egorie des faisceaux pervers sur $\CC^n$ relativement \`a un croisement normal avec une sous-cat\'egorie pleine, not\'ee $\Cc_{n}$,  de la cat\'egorie de repr\'esentations du carquois associ\'e \`a un hypercube. On peut aussi citer \cite{Maiso2}, \cite{BG} et \cite{Nar}.
D'autres m\'ethodes, plus g\'en\'erales mais du m\^eme coup moins ex\-pli\-cites, consistent \`a d\'ecrire la cat\'egorie $\Pc erv_X$ dans le cas o\`u $X$ est un espace de Thom-Mather. Citons par exemple R. MacPherson et K. Vilonen, dans \cite{McV}.  

Les r\'esultats que nous venons de citer ne tiennent pas compte du caract\`ere local des faisceaux pervers. Rappelons, en effet, que les faisceaux pervers sur un espace stratifi\'e forment un champ. Autrement dit conna\^itre un faisceau pervers \`a isomorphisme pr\`es sur un ouvert $U$ revient \`a le conna\^itre sur un recouvrement ouvert. On peut donc s'attendre \`a ce qu'une caract\'erisation locale de la cat\'egorie des faisceaux pervers par des cat\'egories explicites de repr\'esentations de carquois puisse se recoller sur $X$ en une description simple de la cat\'egorie globale des faisceaux pervers. Ceci permettrait, par exemple, de d\'ecrire les faisceaux pervers sur un espace qui est localement un des espaces cit\'es. Dans \cite{Tr1} et \cite{Tr2}, D. Treumann utilise se genre de m\'ethodes pour g\'en\'eraliser le proc\'ed\'e d\'efini dans \cite{McV}. 

Nous nous int\'eressons ici \`a la cat\'egorie $\Pc erv_X$ des faisceaux pervers sur une vari\'et\'e torique lisse $X$ stratifi\'ee par le l'action du tore. Les vari\'et\'es toriques sont des espaces  indiqu\'es car ils sont localement isomorphes \`a des produits de $\CC^k$ et de $\CC^{*l}$ stratifi\'es par le croisement normal.  Il s'agit donc d'utiliser la notion de champs et notamment de champs constructibles pour recoller les descriptions de A. Galligo, M. Granger et Ph. Maisonobe. Nous exposons ici le r\'esultat de cette m\'ethode. Ainsi \`a  partir de la donn\'ee d\'un \'eventail nous d\'efinissons une sous cat\'egorie  pleine de la cat\'egorie de repr\'esentations d'un certain carquois \'equivalente  \`a la cat\'egorie $\Pc erv_X$.

\section{Cat\'egorie des faisceaux pervers sur les vari\'et\'es toriques lisses.}
On note $\Delta$ un \'eventail de $\CC^n$,  $v_{1}, \cdots, v_{k}$ les vecteurs primitifs qui engendrent les c\^ones de dimension un appartenant \`a $\Delta$, et $\Ic$ l'ensemble des parties $I$ de $\{1, \cdots, k\}$ telles que l'enveloppe positive des vecteurs $\{v_i\}_{i \in I}$ soit un c\^one de $\Delta$. On note $X_{I}$ la vari\'et\'e associ\'ee au c\^one $\sigma_{I}$ et $\Pc erv_{I}$ la categorie  des faisceaux pervers sur $X_{I}$ relativement \`a la stratification donn\'ee par l'action du tore.\\
On suppose de plus que $\Delta$ est un \'eventail r\'egulier de $\RR^n$, la vari\'et\'e $X_\Delta$ est donc lisse.
\begin{Def}
Soit $I$ une partie de $\Ic$, on dit que $I$ et $\sigma_I$ sont maximaux s'ils sont maximaux pour l'inclusion dans respectivement $\Ic$ et $\Delta$.
\end{Def}
Pour tout c\^one maximal $\sigma_{I}$ de $\Delta$, on se fixe une base $\Bc_I$ de $\ZZ^n$ qui contient les vecteurs $\{v_{i}\}_{i \in I}$ (l'existence de $\Bc_I$ est assur\'ee par le fait que $\Delta$ est r\'egulier).
\begin{Def}
On note $c_{\Delta}$ le carquois d\'efini par la donn\'ee :
\begin{itemize}
\item[$\bullet$] pour tout c\^one $\sigma_{I}$ de $\Delta$, d'un sommet $s_{I}$ muni de $n-l$ boucles, $l$ \'etant le cardinal maximal des ensembles maximaux de $\Ic$ contenant $I$,
\item[$\bullet$] pour tout couple $(\sigma_{I}, \sigma_{I'})$ de c\^ones de $\Delta$ tels que $\sigma_{I'}$ soit une face de codimension un de $\sigma_{I}$, de deux fl\`eches entre $s_I$ et $s_{I'}$ de sens inverse.
\end{itemize}
\end{Def}
\textbf{Exemples}
\begin{itemize}
\item Le carquois associ\'e \`a la vari\'et\'e torique $\CC^n$ est un hypercube de dimension $n$. De m\^eme le carquois associ\'e \`a $\CC^l\times \CC^{*n-l}$ est un hypercube de dimension $l$ dont les sommets sont munis de $n-l$ boucles.
\item Si $\Delta \in \CC^2$ est la r\'eunion des faces des c\^ones respectivement engendr\'es par les vecteurs  $(e_{1}, e_{2})$ et $-e_{1}-e_{2}$ o\`u $(e_{1},e_{2})$ est la base canonique de $\CC^2$. Le carquois associ\'e est le carquois $(1)$. De m\^eme les carquois associ\'e \`a $\PP_1$ et $\PP_2$ sont respectivement les carquois $(2)$ et $(3)$:
$$
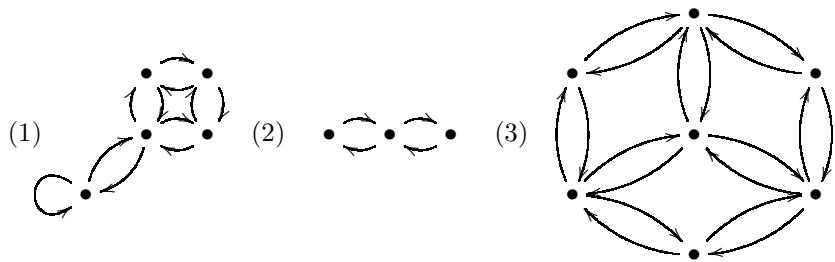
\xymatrix @!0 @C=.8cm @R=.8cm{
&&&&&&&&&&&  \bullet \ar@/^/[dll] \ar@/^/[drr]  \ar@/^/[dd]\\
&& \bullet \ar@/^/[r] \ar@/^/[d] & \bullet \ar@/^/[l] \ar@/^/[d]&&&&&&\bullet \ar@/^/[urr] \ar@/^/[dd]  & &  & & \bullet   \ar@/^/[ull]  \ar@/^/[dd] \\
(1)&& \bullet \ar@/^/[u] \ar@/^/[r] \ar@/^/[dl] & \bullet \ar@/^/[u] \ar@/^/[l]&(2)&\bullet \ar@/^/[r] & \bullet \ar@/^/[r] \ar@/^/[l] &\bullet \ar@/^/[l]&(3) && & \bullet \ar@/^/[uu]  \ar@/^/[dll]  \ar@/^/[drr]\\
&\bullet \ar@/^/[ur] \ar@(ul,dl)[]&&&&&&&& \bullet  \ar@/^/[uu]  \ar@/^/[urr]  \ar@/^/[drr]  & & & & \bullet  \ar@/^/[uu]  \ar@/^/[ull]  \ar@/^/[dll] \\
&&&&&&&&&&&  \bullet  \ar@/^/[ull]  \ar@/^/[urr]
}
$$
Un objet de la cat\'egorie des repr\'esentations du carquois $c_\Delta$ est donc la donn\'ee d'une famille : $$\Big(\{E_I\}_{I \in \Ic}, \{u_{Ip},v_{Ip} \}_{I\cup \{p\} \in \Ic}, \{_iM_{I}\}\Big)$$ o\`u $E_I$ est un espace vectoriel de dimension finie, $u_{Ip}$ et  $v_{Ip}$  sont des applications lin\'eaires  de respectivement $E_I$ dans $E_{Ip}$ et $E_{Ip}$ dans $E_I$ et o\`u $_iM_I$ est un endomorphisme de $E_I$ o\`u $i$ varie de $0$ \`a $n-l$, $l$ \'etant le cardinal maximal des ensembles maximaux de $\Ic$ contenant $I$.
\end{itemize}
\begin{Def}\label{CcDelta}
On note $\Cc_{\Delta}$  la sous-cat\'egorie pleine de la cat\'egorie des repres\'esentations du carquois $c_\Delta$ form\'ee des objets tels que  :
\begin{itemize}
\item[(i)] Pour toute partie $I$ de $\Ic$ l'application lin\'eaire $_iM_I$ est inversible. 
\item[(ii)] Pour tout couple $(I\cup\{p\}, I)$ appartenant \`a $\Ic^2$  l'application lin\'eaire $M_{Ip}$ suivante est inversible :
$$M_{Ip}= v_{Ip} u_{Ip} + Id .$$
\item[(iii)] Pour tout quadruplet $(I, I\cup \{p\}, I\cup \{q\},I \cup \{p, q\})$ de parties appartenant \`a $\Ic$ les applications lin\'eaires $u_{Ip}$, $u_{Iq}$, $u_{Ipq}$, $u_{Iqp}$ et $v_{Ip}$, $v_{Iq}$, $v_{Ipq}$, $v_{Iqp}$ donn\'ees par le diagramme 
$$
\xymatrix @!0 @C=2cm @R=1.4cm{& E_{I\cup p} \ar@/^/[ld]^{v_{Ip}} \ar@/^/[rd]^{u_{Ipq}}&\\
E_{I}\ar@/^/[ur]^{u_{Ip}} \ar@/^/[dr]^{u_{Iq}}& &E_{I \cup \{p, q\}} \ar@/^/[ul]^{v_{Ipq}} \ar@/^/[dl]^{v_{Iqp}}\\
& E_{I \cup q} \ar@/^/[lu]^{v_{Iq}} \ar@/^/[ru]^{u_{Iqp}}
}$$  
 v\'erifient les conditions de commutations suivantes : 
$$u_{Ip} u_{Ipq}= u_{Iq} u_{Iqp}, ~v_{Ipq}v_{Ip}=v_{Iqp} v_{Iq}, ~ v_{Ipq}u_{Iqp}=u_{Ip} v_{Iq}.$$
\item[(iv)] Pour toute partie $J$ de $\Ic$, pour tout \'el\'ement $p$   de $\{1,\cdots,k\}\backslash J$ tel que $J\cup\{p\} \in \Ic$ et pour toute  partie maximale  $K$  de $\Ic$ de cardinal maximal contenant $J$ mais ne contenant pas $p$,  on note $(v_1, \cdots, v_n)$ la base $\Bc_{K}$ fix\'ee plus haut et on suppose que, pour tout $i \in \{1, \cdots,j\}$ $v_i $ appartient \`a $\sigma_J$, pour tout $i \in \{j+1, \cdots,l\}$ $v_i $ appartient \`a $\sigma_K$. Enfin on note $(p_1, \cdots, p_n)$ les coordonn\'ees de $p$ dans $\Bc_K$. Alors l'application lin\'eaire $M_{Jp}$ s'\'ecrit comme suit :
$$M_{Jp}= M_{Jv_{j+1}}^{p_{j+1}} \cdots  M_{Jv_{l}}^{p_{l}}(_1M_J)^{p_{l+1}}\cdots (_{n-l}M_J)^{p_{n}}$$
\end{itemize}
\end{Def}
Pour mieux comprendre cette d\'efinition donnons quelques exemples : 
\begin{itemize}
\item[$\bullet$]
Quand la vari\'et\'e torique est $\CC^n$ on retrouve la cat\'egorie $\Cc_{\CC^n}$ d\'efinie par G. Galligo, M. Granger et Ph. Maisonobe dans \cite{GGM}.
\item[$\bullet$] Soit $(e_1, e_2)$ la base canonique de $\ZZ^2$, $e_3=-e_1-e_2$ et $\Delta$ l'\'eventail form\'e des faces des c\^ones engendr\'es par $\{e_1,e_2\}$, $\{e_1,e_3\}$ et $\{e_2, e_3\}$. La vari\'et\'e torique associ\'ee \`a $\Delta$ est $\PP^2$.
La cat\'egorie $\Cc_{\PP^2}$ est la cat\'egorie form\'ee des objets : 
$$
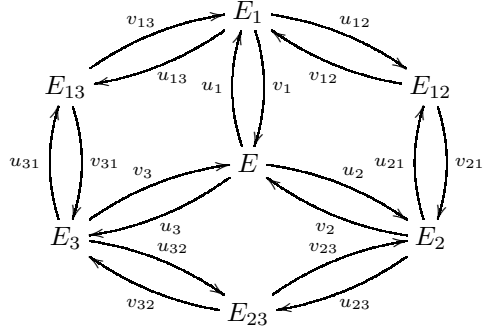
\xymatrix @!0 @C=1.2cm @R=1cm{    & & E_1 \ar@/^/[dll]^{u_{13}}  \ar@/^/[drr]^{u_{12}}  \ar@/^/[dd]^{v_{1}}  \\
                           E_{13}  \ar@/^/[urr]^{v_{13}}  \ar@/^/[dd]^{v_{31}}  & &  & & E_{12}   \ar@/^/[ull]^{v_{12}}  \ar@/^/[dd]^{v_{21}} \\
                           & & E \ar@/^/[uu]^{u_1}  \ar@/^/[dll]^{u_{3}}  \ar@/^/[drr]^{u_{2}}   \\
                           E_3  \ar@/^/[uu]^{u_{31}}  \ar@/^/[urr]^{v_{3}}  \ar@/^/[drr]^{u_{32}}   & & & & E_2  \ar@/^/[uu]^{u_{21}}  \ar@/^/[ull]^{v_{2}}  \ar@/^/[dll]^{u_{23}}   \\
                           & & E_{23}  \ar@/^/[ull]^{v_{32}}  \ar@/^/[urr]^{v_{23}}  \\   }$$  
tels que, pour $(i,j)\in \{1,2\}$ avec $i\neq j$ :
\begin{itemize}
\item[-]
$u_{i} u_{ij}= u_{j} u_{ji}, ~v_{ij}v_{i}=v_{ji} v_{j}, ~ v_{ji}u_{ij}=u_{j} v_{i}$
\item[-] $M_{\emptyset i}= v_i \circ u_i+ Id$ et $M_{ij}= v_{ij} \circ u_{ij} + Id$ soient inversibles,
\end{itemize}
Explicitons maintenant la condition $(iv)$ :
\begin{itemize}
\item[-] Consid\'erons  $M_{\emptyset 3}$, ici $J=\emptyset$, $p=\{e_{3}\}$ et $ K=\Bc_K=(e_1,e_2)$ ainsi on a :
$M_{\emptyset3} =M_{\emptyset1}^{-1} M_{\emptyset2}^{-1},$
\item[-]  Int\'eressons nous maintenant \`a l'endomorphisme $M_{13}$, cette fois $J=\{e_1\}$, $p=\{e_{3}\}$ et $ K= \Bc_K=(e_1,e_2)$ on a ainsi :
$M_{13}=M_{12}^{-1}$ 
De m\^eme, pour tout $(i,j,k) \in \{1,2,3\}^3$, tous les trois diff\'erents on a :
$M_{ij}= M_{ik}^{-1}$
\end{itemize}   
\end{itemize}                                               
\begin{thm}\label{torsection}
Les cat\'egories $\Cc_{\Delta}$ et $\Pc erv_{\Delta}$ sont \'equivalentes.
\end{thm}
~\\
Soit $\Delta$ un \'eventail r\'egulier et $X$ la vari\'et\'e torique associ\'ee \`a cet \'eventail.
L'id\'ee de la preuve est de d\'efinir un champ $\CCC_\Delta$ de cat\'egories de repr\'esentations de carquois \'equivalent au champ des faisceaux pervers sur $X$ relativement \`a la stratification de l'action du tore. En effet la cat\'egorie des sections globales de $\CCC_X$ est alors \'equivalente \`a la cat\'egorie $\Pc erv_\Delta$. Les champs pouvant \^etre d\'efinis localement il suffit de d\'efinir un tel champ pour toute vari\'et\'e torique lisse affine, puis de les recoller.

Notons que les vari\'et\'e toriques lisses affines sont des produits de $\CC^{*}$ et de $\CC$ et que la stratification de l'action du tore n'est autre que la stratification produit du croisement normal et des $\CC^{*}$. Ainsi dans un premier temps nous transformons l'\'equivalence de cat\'egories de  A. Galligo, M. Granger et Ph. Maisonobe en une \'equivalence de champs sur $\CC^n$. Nous nous appuyons pour cela sur une description combinatoire des champs constructibles. 

Puis pour recoller ces champs, et surtout pour expliciter les sections globales du champ recoll\'e nous explicitons, \`a partir de la donn\'ee d'un \'eventail, les morphismes de recollement des vari\'et\'es toriques.

Les conditions $(i)$, $(ii)$ et $(iii)$ de la cat\'egorie $\Cc_\Delta$ sont issues des \'equivalences d\'efinies par A. Galligo, M. Granger et Ph. Maisonobe. La condition $(iv)$ est elle donn\'ee par les recollements toriques. 

\section{Remerciements}
Cette note est issue de ma th\`ese \cite{these}. Je remercie mon directeur Ph. Maisonobe pour sa disponibilit\'e et I. Waschkies pour les nombreuses discussions. Je remercie aussi P. Abbrugiati pour un certain lemme d'alg\`ebre lin\'eaire.

\bibliographystyle{plain}
\bibliography{biblio}

\begin{thebibliography}{1}

\bibitem{BBD}
A.~Beilinson, J.~Bernstein, and P.~Deligne.
\newblock Faisceaux pervers.
\newblock {\em Ast\'erisque}, 100, 1982.

\bibitem{BG}
T.~Braden and M.~Grinberg.
\newblock Perverse sheaves on rank stratifications.
\newblock {\em Duke math. jour.}, 96:317--361, 1999.

\bibitem{these}
D.~Dupont.
\newblock {\em Exemples de classifications du champ des faisceaux pervers}.
\newblock PhD thesis, 2008.

\bibitem{GGM}
A.~Galligo, M.~Granger, and Ph. Maisonobe.
\newblock $\mathcal{D}$-modules et faisceaux pervers dont le support singulier
  est un croisement normal.
\newblock {\em Ann. Inst. Fourier (Grenoble)}, 1985.

\bibitem{Nar}
F.~Gudiel~Rodr\`iguez and L.~Narv\`aez~Macarro.
\newblock Explicit models for perverse sheaves.proceedings of the international
  conference on algebraic geometry and singularities (spanish) (sevilla, 2001).
\newblock {\em Rev. Mat. Iberoamericana}, 19:425--454, 2003.

\bibitem{McV}
R.~MacPherson and K.~Vilonen.
\newblock Elementary construction of perverse sheaves.
\newblock {\em Invent. Math.}, 84:403--435, 1986.

\bibitem{Maiso2}
Ph. Maisonobe.
\newblock Faisceaux pervers dont le support singulier est une courbe plane.
\newblock {\em Compos. Math.}, 62:215--261, 1987.

\bibitem{Tr2}
D.~Treumann.
\newblock Stacks similar to the stack of perverse sheaves.
\newblock {\em http://arxiv.org/pdf/0801.3016}, 2008.

\bibitem{Tr1}
D.~Treumann.
\newblock Exit paths and constructible stacks.
\newblock {\em Compos. Math.}, 145:1504--1532, 2009.

\end{thebibliography}

\end{document}